\documentclass[12pt]{article}
\usepackage{amssymb,amsmath,amsfonts}
\usepackage{amsthm}
\usepackage{color}

\begin{document}

\centerline{\bf\large On notions of Q-independence and Q-identical
distributiveness}
\vspace{0.5cm}

\centerline{\bf Alexander Il'inskii}

\begin{quote}
In a recent article A.M.~Kagan and G.J.~Sz\'ekely introduced a notion of $Q$-independent and $Q$-identical distributed random variables. We give a complete description of polynomials
which appear in these definitions.
\end{quote}

\noindent {\small \textbf{Keywords:} $Q$-independence, $Q$-identical
distributiveness, random variable, characteristic function
\vspace{0.2cm}}

\vspace{0.2cm}

{\bf 1. Introduction}
\vspace{0.2cm}

In the article \cite{KS} A.M.~Kagan and G.J.~Sz\'ekely  introduced a notion of $Q$-independent and  $Q$-identical  distributed
random variables.
\vspace{0.2cm}

{\bf Definition 1.} Let $(X_1,\ldots,X_n)$ be a random vector  with the characteristic function $\varphi(t_1,\ldots,t_n)$. Denote by
$\varphi_1(t),\ldots,\varphi_n(t)$ the characteristic functions of the random variables $X_1,\ldots,X_n$ respectively.
The random variables $X_1,\ldots,X_n$ are said to be {\it $Q$-independent} if for all $t_1,\ldots,t_n\in\mathbb{R}$
the  condition
\begin{equation}\label{1}
\varphi(t_1,\ldots,t_n)=\varphi_1(t_1)\cdots\varphi_n(t_n)
\exp(q(t_1,\ldots,t_n)),
\end{equation}
where $q(t_1,\ldots,t_n)$ is
a polynomial such that $q(0,\ldots,0)=0$, holds.
\vspace{0.2cm}

{\bf Definition 2.}  Random variables $X$ and $Y$ with the characteristic functions $\varphi_X(t)$ and $\varphi_Y(t)$ respectively are said to be
{\it $Q$-identically distributed} if for all $t\in\mathbb{R}$ the condition
\begin{equation}\label{2}
\varphi_X(t)=\varphi_Y(t)\exp(q(t)),
\end{equation}
where  $q(t)$ is a polynomial such that $q(0)=0$,
holds.

A.M.~Kagan and G.J.~Sz\'ekely \cite{KS} proved that a series of classical theorems (Cram\'er \cite{C}, Marcinkiewicz \cite{Mar}, Skitovich--Darmois  \cite{S}, \cite{D}, and Vershik \cite{V}) holds true if the conditions of independence and identical distributiveness are changed for the conditions of
$Q$-independence and $Q$-identical distributiveness. The paper \cite{KS} stimulated the appearance of the papers \cite{F}, \cite{M}, \cite{R1}, \cite{R2}.
B.L.S. Prakasa Rao \cite{R1}, \cite{R2} considered some generalizations of the Kotlarski \cite{K}, \cite{RR} theorem and
similar results for Q-independent random variables. Following A.M.~Kagan and G.J.~ Sz\'ekely \cite{KS} G.M.~Feldman in \cite{F} introduced a
notion of Q-independence for random variables taking values in a locally compact Abelian group.
He proved in \cite{F} that if we consider Q-independence instead
of independence, then the group analogue
of the Cram\'er theorem (\cite{F1}, \cite{F2}), and some group analogues of the the
Skitovich--Darmois (\cite{F3}) and
Heyde theorems (\cite{F4}) hold true for the same classes of groups.
M.V.~Myronyuk in \cite{M} continues research of characterization theorems for Q-independent random
variables with values in a locally compact Abelian group.

Taking into account Definitions~1 and~2 the following question naturally arises: what kind of polynomials can be involved in the definitions of $Q$-independence and $Q$-identical distributiveness in formulas \eqref{1} and \eqref{2}? In \cite{KS} authors describe a wide class of polynomials for Definition~2. They also give  an important example for Definition~1 with polynomial $q(t_1,t_2)=t_1t_2$.
The aim of this paper is to give a complete description of the class  of polynomials that can be involved in definitions of $Q$-independence and $Q$-identical distributiveness.
\vspace{0.2cm}

{\bf 2. $Q$-independence}
\vspace{0.2cm}

{\bf Theorem 1.} {\it Let
$$
q(t_1,\dots,t_n)=\sum\limits_{k_1,\dots,k_n\geq 0}a_{k_1,\dots,k_n}t_1^{k_1}\ldots t_n^{k_n}\quad(a_{0,\dots,0}=0)
$$
be a polynomial of $n$ variables. Then there exists a random vector $(X_1,\ldots,X_n)$ such that for all $t_1,\ldots,t_n\in\mathbb{R}$   equality $(\ref{1})$
holds if and only if the following conditions are satisfied$:$

$(i)$ The polynomial $q(t_1,\dots,t_n)$ does not contain monomials which depend only on one variable,

$(ii)$ The coefficient $a_{k_1,\ldots,k_n}$ of the polynomial $q(t_1,\dots,t_n)$ is real if $k_1+\ldots+k_n$ is even, and $a_{k_1,\dots,k_n}$ is pure imaginary
if $k_1+\ldots+k_n$ is odd.}
\vspace{0.1cm}

We note that our proof of Theorem~1 does not depend on $n$.
For  brevity we give the proof for $n=2$.

{\bf Proof. Necessity.} Put $\|t\|=(t_1^2+t_2^2)^{1/2}$. Let $\varphi(t_1,t_2)$ be the characteristic function of the random vector $(X_1,X_2)$.  Denote by $\varphi_1(t)$ and $\varphi_2(t)$  the characteristic functions of the random variables  $X_1$ and $X_2$ respectively.  Suppose that for all $t_1,t_2\in \mathbb{R}$ we have
\begin{equation}\label{4}
\varphi(t_1,t_2)=\varphi_1(t_1)\varphi_2(t_2)\exp(q(t_1,t_2))\,.
\end{equation}
It is obvious that $\varphi(0,t_2)=\varphi_2(t_2)$, $t_2\in \mathbb{R}$.  It follows from \eqref{4} for $t_1=0$ that if $\delta>0$ is small enough, and
$t_2\in \mathbb{R}$ is such that $\|t\|<\delta$  we have $\exp(q(0,t_2))=1$.
We choose $\delta$  so small that for $\|t\|<\delta$ all characteristic functions in \eqref{4} do not vanish.
Taking into account that $q(0,0)=0$, we obtain that $q(0,t_2)=0$.
Therefore $a_{0,k_2}=0$ $(k_2\ge 1)$.
Reasoning similarly we get $a_{k_1,0}=0$ $(k_1\ge 1)$.
This proves  that condition $(i)$ of Theorem~1 holds.

Condition $(ii)$ of Theorem~1
follows from the Hermitian property of the characteristic functions. Indeed, we have $\varphi(-t_1,-t_2)=\overline{\varphi(t_1,t_2)}$,
and  $\varphi_j(-t)=\overline{\varphi_j(t)}$, $j=1,2$, for all $t_1,t_2\in \mathbb{R}$. This implies from \eqref{4}
that $\exp(q(-t_1,-t_2))=\exp(\overline{q(t_1,t_2)})$ for $\|t\|<\delta$, where $\delta>0$ is sufficiently small.
From $q(0,0)=0$ it follows that
$q(-t_1,-t_2)=\overline{q(t_1,t_2)}$.
Therefore we have for $\|t\|<\delta$
$$
\sum\limits_{k_1,k_2\geq 1}(-1)^{k_1+k_2}a_{k_1,k_2}t_1^{k_1}t_2^{k_2}=
\sum\limits_{k_1,k_2\geq 1}\overline{a}_{k_1,k_2}t_1^{k_1}t_2^{k_2}\,.
$$
This means that $(-1)^{k_1+k_2}a_{k_1,k_2}=\overline{a}_{k_1,k_2}$
for all $k_1,k_2\geq 1$.
This is equivalent to condition~$(ii)$ of Theorem~1.
\vspace{0.2cm}

{\bf Sufficiency.} The proof of sufficiency is based on the following theorem of A.A.~Goldberg \cite{G}.
\vspace{0.2cm}

{\bf Theorem A.}
{\it There exists an entire function $f(z)$ such that for some positive constants $C_1$, $C_2$, $C_3$, all $k\in\mathbb{N}$ and $x\in\mathbb{R}$ the conditions$:$

$1)$ $0<f(x)<C_1\exp(-C_2\log^2(|x|+1))$,

$2)$ $|f^{(k)}(x)|\leq C_3^kf(x)$\\
hold.}
\vspace{0.2cm}

Obviously, the probability density $f(x)/\int_{\mathbb{R}}f(v)dv$
also satisfies conditions~$1)$ and~$2)$ of Theorem~A (with another constant $C_1$).
Therefore, without loss of generality we can assume  that the function $f(x)$ itself is a probability density.

Let $0<\varepsilon<1$.
Consider the new probability density
\begin{equation}\label{5}
p_{\varepsilon}(x):=\varepsilon f(\varepsilon x)\,.
\end{equation}
Denote  by
\begin{equation}\label{55}
\varphi_{\varepsilon}(t):=
\int\limits_{\mathbb{R}}e^{itx}p_{\varepsilon}(x)dx\
\end{equation}
the characteristic function of $p_{\varepsilon}(x)$.

We will show that if  $\varepsilon$ is small enough, then the function
\begin{equation}\label{6}
e^{q(t_1,t_2)}\varphi_{\varepsilon}(t_1)\varphi_{\varepsilon}(t_2)
\end{equation}
is the characteristic function $\varphi(t_1,t_2)$ of a random vector $(X_1,X_2)$. This proves Theorem~1. Indeed,
if we take $t_1=0$ in the equality
\begin{equation}\label{7}
\varphi(t_1,t_2)=e^{q(t_1,t_2)}\varphi_{\varepsilon}(t_1)\varphi_{\varepsilon}(t_2)\,,
\end{equation}
then using condition $(i)$ of Theorem~1 we obtain
$$
\varphi_{\varepsilon}(t_2)=\varphi(0,t_2)=\varphi_{2}(t_2)\,.
$$
Reasoning similarly we get
$$
\varphi_{\varepsilon}(t_1)=\varphi(t_1,0)=\varphi_{1}(t_1)\,.
$$
Hence, $(\ref{7})$ coincides with $(\ref{4})$.

We follow now the paper \cite{G}. In sequel we use notation $\displaystyle D_x^m:=\frac{d^m}{dx^m}$.
Integrating by parts we obtain from $(\ref{55})$ that for every $m\in\mathbb{N}$ the following equality
\begin{equation}\label{8}
t^m\varphi_{\varepsilon}(t)=
t^m\int\limits_{\mathbb{R}}e^{itx}p_{\varepsilon}(x)dx
=i^m\int\limits_{\mathbb{R}}e^{itx}D_x^mp_{\varepsilon}(x)dx
\end{equation}
holds.
It follows from  $(\ref{5})$ and statement $2)$ of Theorem~A that for all $m\in\mathbb{N}$ and $x\in\mathbb{R}$ the inequality
\begin{equation}\label{9}
\left|D_x^mp_{\varepsilon}(x)\right|\leq\varepsilon^mC_3^mp_{\varepsilon}(x)
\end{equation}
is valid, where $C_3$ is the constant from statement $2)$ of Theorem~A.
We want to prove that the function \eqref{6} is the characteristic function of a probability density $r(x_1,x_2)$, i.e.
\begin{equation}\label{10}
e^{q(t_1,t_2)}\varphi_{\varepsilon}(t_1)\varphi_{\varepsilon}(t_2)=
\iint\limits_{\mathbb{R}^2}e^{i(t_1x_1+t_2x_2)}r(x_1,x_2)dx_1dx_2\,.
\end{equation}
We have
\begin{equation}\label{11}
e^{q(t_1,t_2)}\varphi_{\varepsilon}(t_1)\varphi_{\varepsilon}(t_2)=
\iint\limits_{\mathbb{R}^2}e^{i(t_1x_1+t_2x_2)}p_{\varepsilon}(x_1)p_{\varepsilon}(x_2)dx_1dx_2+
\sum\limits_{n=1}^{\infty}\frac{q^n(t_1,t_2)}{n!}\varphi_{\varepsilon}(t_1)\varphi_{\varepsilon}(t_2)\,.
\end{equation}
Consider the first term ($n=1$) in series $(\ref{11})$. Using $(\ref{8})$ we obtain
\begin{equation*}
\begin{split}
q(t_1,t_2)\varphi_{\varepsilon}(t_1)\varphi_{\varepsilon}(t_2)&=
\sum\limits_{k_1,k_2\geq 1}a_{k_1,k_2}t_1^{k_1}
\int\limits_\mathbb{R} e^{it_1x_1}p_{\varepsilon}(x_1)dx_1\cdot
t_2^{k_2}\int\limits_\mathbb{R} e^{it_2x_2}p_{\varepsilon}(x_2)dx_2\\
&=\sum\limits_{k_1,k_2\geq 1}a_{k_1,k_2}i^{k_1+k_2}
\int\limits_{\mathbb{R}}e^{it_1x_1}D_{x_1}^{k_1}p_{\varepsilon}(x_1)dx_1\cdot
\int\limits_{\mathbb{R}}e^{it_2x_2}D_{x_2}^{k_2}p_{\varepsilon}(x_2)dx_2\\
&=\iint\limits_{\mathbb{R}^2}e^{i(t_1x_1+t_2x_2)}\left(\sum\limits_{k_1,k_2\geq 1}a_{k_1,k_2}i^{k_1+k_2}
D_{x_1}^{k_1}p_{\varepsilon}(x_1)D_{x_2}^{k_2}p_{\varepsilon}(x_2)\right)dx_1dx_2\\
&=\iint\limits_{\mathbb{R}^2}e^{i(t_1x_1+t_2x_2)}S_1(x_1,x_2)dx_1dx_2,
\end{split}
\end{equation*}
where
\begin{equation}\label{11_1}
S_1(x_1,x_2)=\sum\limits_{k_1,k_2\geq 1}a_{k_1,k_2}i^{k_1+k_2}
D_{x_1}^{k_1}p_{\varepsilon}(x_1)D_{x_2}^{k_2}p_{\varepsilon}(x_2).
\end{equation}
It follows from condition $(ii)$ of Theorem~1 that
$a_{k_1,k_2}i^{k_1+k_2}\in\mathbb{R}$ for all $k_1,k_2$. This implies that
$S_1(x_1,x_2)$ is a real valued function.
Let us estimate $|S_1(x_1,x_2)|$.
In what follows we denote by $s$ the number of monomials in the polynomial $q$, and put
$$
d:=\deg q=\max\{k_1+k_2:a_{k_1,k_2}\not=0\}\,,\quad a:=\max\{|a_{k_1,k_2}|:k_1,k_2\}\,.
$$
Note that we may suppose that  $C_3>1$ in inequality \eqref{9}. Using inequality \eqref{9} and taking into account that $0<\varepsilon<1$
we obtain from $(\ref{11_1})$
\begin{equation*}
\begin{split}
|S_1(x_1,x_2)|&\leq
\sum\limits_{k_1,k_2\geq 1}|a_{k_1,k_2}|
\left|D_{x_1}^{k_1}p_{\varepsilon}(x_1)\right|\cdot\left|D_{x_2}^{k_2}p_{\varepsilon}(x_2)\right|\\
&\leq\sum\limits_{k_1,k_2\geq 1}a\varepsilon^{k_1+k_2}C_3^{k_1+k_2}p_{\varepsilon}(x_1)p_{\varepsilon}(x_2)\leq
sa\varepsilon^2C_3^dp_{\varepsilon}(x_1)p_{\varepsilon}(x_2)\,.
\end{split}
\end{equation*}

Thus,
$$
\left|S_1(x_1,x_2)\right|\leq sa\varepsilon^2C_3^dp_{\varepsilon}(x_1)p_{\varepsilon}(x_2)\,.
$$
Reasoning similarly we get
$$
q^n(t_1,t_2)\varphi_{\varepsilon}(t_1)\varphi_{\varepsilon}(t_2)=
\iint\limits_{\mathbb{R}^2}e^{i(t_1x_1+t_2x_2)}S_n(x_1,x_2)dx_1dx_2\
$$
for all $n\in\mathbb{N}$, where
$S_n(x_1,x_2)$ is a real valued function, and moreover, for all $(x_1,x_2)\in\mathbb{R}^2$ we have the inequality
\begin{equation}\label{14}
\left|S_n(x_1,x_2)\right|\leq\left(sa\varepsilon^2C_3^d\right)^np_{\varepsilon}(x_1)p_{\varepsilon}(x_2)\,.
\end{equation}
Therefore \eqref{10} is valid, where
$$
r(x_1,x_2):=p_{\varepsilon}(x_1)p_{\varepsilon}(x_2)+\sum\limits_{n=1}^{\infty}\frac{S_n(x_1,x_2)}{n!}
$$
is a real valued function. Using \eqref{14} and assuming that $\varepsilon$ is small enough we get
\begin{equation*}
\begin{split}
r(x_1,x_2)&=p_{\varepsilon}(x_1)p_{\varepsilon}(x_2)+
\sum\limits_{n=1}^{\infty}\frac{S_n(x_1,x_2)}{n!}
\geq
p_{\varepsilon}(x_1)p_{\varepsilon}(x_2)-\sum\limits_{n=1}^{\infty}\frac{|S_n(x_1,x_2)|}{n!}\\
&\geq p_{\varepsilon}(x_1)p_{\varepsilon}(x_2)\Big(1-\sum\limits_{n=1}^{\infty}\frac{1}{n!}\left(\varepsilon^2saC_3^d\right)^n\Big)\\
&=p_{\varepsilon}(x_1)p_{\varepsilon}(x_2)\big(1-\left\{\exp(\varepsilon^2saC_3^d)-1\right\}\big)>0\,.
\end{split}
\end{equation*}
Thus, the function $r(x_1,x_2)$
is  positive on  $\mathbb{R}^{2}$.  The left-hand side of \eqref{10} is equal to~$1$ for $t_1=t_2=0$. Therefore, the function $r(x_1,x_2)$ is a probability density.
Theorem~1 is completely proved.
\vspace{0.2cm}

{\bf 3. $Q$-identical distributiveness}
\vspace{0.2cm}

{\bf Theorem 2.} {\it Let
$$
q(t)=\sum\limits_{k\geq 1}a_{k}t^k
$$
be a polynomial. Then there exist random variables $X$ and $Y$ such that
for all $t\in\mathbb{R}$ equality \eqref{2} holds if and only if the following conditions are valid$:$

$(ii')$ $a_k\in\mathbb{R}$ if $k$ is even and $a_k\in i\mathbb{R}$ if $k$ is odd.}
\vspace{0.2cm}

The proof of Theorem~2 is completely analogous to the proof of Theorem~1. Arguing as in proving of necessity of Theorem~1 we obtain condition $(ii')$. For the proof of sufficiency in Theorem~2 we need to show that the function $e^{q(t)}\varphi_{\varepsilon}(t)$ (see \eqref{55}) is a characteristic function. We omit this reasoning.
\vspace{0.2cm}

The author expresses sincere gratitude to G.M.~Feldman for suggesting this problem.

{\small
\noindent
V.N. Karazin Kharkiv National University\\
\noindent
Department of Mathematics and Informatics\\
\noindent
61022, Kharkiv, Ukraine\\}
\noindent
{\it e-mail: iljinskii@univer.kharkov.ua}\\

\end{document}